\documentclass[reqno,12pt]{amsart}
\usepackage{setspace,tikz,xcolor,mathrsfs,listings,multicol,amssymb,amsfonts}
\usepackage{pgfplots}
\usepackage{rotating}
\usepackage[vcentermath]{youngtab}
\usepackage{tkz-base}
\usepackage{tkz-euclide}
\usepackage[margin=1in,includefoot,footskip=30pt]{geometry}
\usepackage{enumerate}
\usepackage{booktabs,subcaption}
\usepackage{gensymb}
\usepackage{ytableau}
\usepackage[all,cmtip]{xy}
\usetikzlibrary{arrows,matrix,math}
\tikzset{tab/.style={matrix of math nodes,column sep=-.35, row sep=-.35,text height=7pt,text width=7pt,align=center,inner sep=2,font=\footnotesize}}
\newif\iftikz
\tikztrue

\newif\iflimitshapes
\limitshapestrue  % comment out to compile faster

\usepackage[colorlinks=true, pdfstartview=FitV, linkcolor=blue, citecolor=blue, urlcolor=blue]{hyperref}

\newcommand{\gl}{\mathfrak{gl}}

\newcommand{\NN}{\mathbb{N}}
\newcommand{\bigO}{\mathcal{O}} % big O
\DeclareMathOperator{\Ai}{Ai} % Airy function
\DeclareMathOperator{\Pf}{Pf}
\newcommand{\mcC}{\mathcal{C}}

\newcommand{\II}{\mathbb{I}}
\newcommand{\PP}{\mathbb{P}}
\newcommand{\ZZ}{\mathbb{Z}}

\newcommand{\RR}{\mathbb{R}}
\newcommand{\CC}{\mathbb{C}}

\newcommand{\la}{\lambda}

\theoremstyle{plain}
\newtheorem{thm}{Theorem}[section]

\newtheorem{prop}[thm]{Proposition}

\theoremstyle{definition}

\numberwithin{equation}{section}

\begin{document}  
  
\author[D.~Betea]{Dan Betea}
\address[D.~Betea]{Université d’Angers, CNRS, LAREMA, SFR MATHSTIC, Angers, F-49045, France}
\email{dan.betea@gmail.com}
\urladdr{https://sites.google.com/view/danbetea}

\author[A.~Nazarov]{Anton Nazarov}
\address[A.~Nazarov]{Department of High Energy and Elementary Particle Physics, St.\ Petersburg State University, University Embankment, 7/9, St.\ Petersburg, Russia, 199034 and Beijing Institute of Mathematical Sciences and Applications (BIMSA),
Bejing 101408, People’s Republic of China}
\email{antonnaz@gmail.com}
\urladdr{http://hep.spbu.ru/index.php/en/1-nazarov}

\author[P.~Nikitin]{Pavel Nikitin}
\address[P.~Nikitin]{%Laboratory of Representation Theory and Dynamical Systems, St.\ Petersburg Department of Steklov Mathematical Institute of Russian Academy of Sciences, 27 Fontanka, St.\ Petersburg, Russia, 191023.
  Beijing Institute of Mathematical Sciences and Applications (BIMSA),
Bejing 101408, People’s Republic of China}
\email{pnikitin0103@yahoo.co.uk}

%\author[T.~Scrimshaw]{Travis Scrimshaw}
%\address[T.~Scrimshaw]{Department of Mathematics, Hokkaido University, 5 Ch\=ome Kita 8 J\=onishi, Kita Ward, Sapporo, Hokkaid\=o 060-0808}
%\email{tcscrims@gmail.com}
%\urladdr{https://tscrim.github.io/}

%\title{Non-symmetric analogue of Schur {\color{red}??Demazure} measure and symmetric last passage percolation}

\title{Last passage percolation in lower triangular domain}

\begin{abstract}
  Last passage percolation (LPP) in an $n\times n$ lower triangular domain has nice connections with various generalizations of Schur measures. LPP along an anti-diagonal, from $(1,n)$ to $(n,1)$, gives a distribution of a highest column of a random composition with respect to a Demazure measure (a non-symmetric analog of a Schur measure). LPP along a main diagonal, from $(1,1)$ to $(n,n)$, is distributed as a marginal of a Pfaffian Schur process. In the first case we show that the asymptotics for the constant specialization is governed by the GOE Tracy--Widom distribution, in the second case~--- by the GSE Tracy--Widom distribution. In the latter case we were also able to study the truncated lower triangular case, obtaining an interesting generalization of the GSE Tracy--Widom distribution.
%  Classical Cauchy identity is a manifestation of $(\gl_{n},\gl_{m})$
%  Howe duality. Restriction to Borel subalgebras produces a
%  non-symmetric analogue of the Cauchy identity written in terms of
%Demazure characters and Demazure atoms enumerated by compositions
%$\mu\in \ZZ^n_{\ge0}$. This defines a Demazure measure as a non-symmetric analogue of a Schur measure. We consider a {\color{red}restriction} of the measure to the staircase and truncated staircase shapes and demonstrate that the asymptotic behavior of the top entry of a random composition is governed by GSE Tracy-Widom distribution and its generalization.
\end{abstract}

\maketitle

\section{Introduction}

Consider the following last passage percolation (LPP) model. Take two sets of non-negative parameters $\{x_{1},\dotsc,x_{n}\}$, $\{y_{1},\dotsc,y_{m}\}$, $m, n\in\NN$, such that $x_{i}y_{j}<1$ for all $i,j$. Let $\{w_{i,j}\}$, $1\le i\le n$, $1\le j\le m$ be a family of independent geometric random variables,
\begin{equation}\label{eq:square-LPP-def}
    \PP(w_{i,j} = k) = (1-x_i y_j) (x_i y_j)^k, \quad k\in \ZZ_{\ge 0},
\end{equation}
filling an $n\times m$ matrix $A$. Consider down-right paths $\pi=\{\pi(p)\}$ from the position $(1,1)$ to $(n,m)$, i.e. $\pi(p) - \pi(p-1)\in\{(0,1),(1,0)\}$. The LPP time is then defined as
\begin{equation*}
L_{(1,1)\to(n,m)} = \max_{\pi} \sum_{(r,t)\in\pi}w_{r,t},  
\end{equation*}
with the maximum taken over all the corresponding down-right paths $\pi$. The natural question is to understand the asymptotic behavior of the LPP time as $n, m$ tends to infinity.

Applying the classical  Robinson--Schensted--Knuth (RSK) algorithm to the random matrix $A$ we obtain a pair of semistandard Young tableaux of shape $\lambda$  \cite{fulton1997young,schensted1961longest,knuth1970permutations}. Random Young diagram $\lambda$ is then distributed with respect to the Schur measure~\cite{Okounkov01}:
\begin{equation*}
  \PP^{\text{Schur}}(\lambda)=s_\la(x_{1},\dotsc, x_{n}) s_\la(y_{1},\dotsc,y_{m}) \prod_{i=1}^n \prod_{j=1}^m (1-x_i y_j) .
\end{equation*}
 Here $s_{\lambda}(x_{1},\dotsc, x_{n})$ are Schur polynomials that can be written as sums over semistandard Young tableaux
\begin{equation*}
s_{\lambda}(x_{1},\dotsc, x_{n})=\sum_{T\in SSYT(\lambda|n)}\prod_{i=1}^{n}x_{i}^{T_{i}},  
\end{equation*}
and the normalization constant is obtained from the classical Cauchy identity,  
\begin{equation*}
  \label{eq:cauchy-identity}
\sum_{\lambda} s_\la(x_{1},\dotsc, x_{n}) s_\la(y_{1},\dotsc,y_{m}) = \prod_{i=1}^n \prod_{j=1}^m \frac1{1-x_i y_j}.  
\end{equation*}
The length of the first row $\lambda_{1}$ is equal to the LPP time. The  question of its asymptotics for a certain choice of parameters was initially posed by Ulam as the question about the length of a longest increasing subsequence in a random permutation. It was famously solved by \cite{vershik1977asymptotics,logan1977variational} with further elaboration by \cite{Johansson00, BDJ99}. 

The classical Cauchy identity is a manifestation of $(\gl_{n},\gl_{m})$ Howe duality \cite{howe1989remarks}. One of the main motivations for the present paper was the study of its non-symmetric analogue, implicitly introduced by van der Kallen in~\cite{vdK89}. Restriction to Borel subalgebras produces a non-symmetric analogue of the Cauchy identity written in terms of Demazure characters and Demazure atoms enumerated by compositions $\mu\in \ZZ^n_{\ge0}$. The corresponding analogue of RSK algorithm was given by S.~Mason in \cite{Mason2006}. From an $n\times m$ matrix $A$ with coefficients from $\ZZ_{\geq 0}$ the algorithm produces a pair of semi-skyline augmented filings (SSAF), which play the role of semistandard Young tableaux in this case. The shapes of SSAFs are compositions
$\mu\in \ZZ^n_{\ge0}$ and we get an analogue of the Cauchy identity
\begin{equation*}
  \sum_{\mu\in\ZZ^{n}_{\geq 0}} \bar{\kappa}^{\mu}(x_{1},\dotsc,x_{n}) \kappa_{\sigma_0\la(\mu)}(y_{1},\dotsc,y_{m})=\prod_{i=1}^{n}\prod_{j=1}^{m}(1-x_{i}y_{j}), 
\end{equation*}
where $\kappa_\mu$ and $\bar{\kappa}^\mu$ are Demazure characters
and opposite Demazure atoms correspondingly, $\la(\mu)$ is a  permutation of the columns of $\mu$ in non-increasing order, and $\sigma_0$ is the longest element of the symmetric group $S_n$. Note that in this case $\kappa_{\sigma_0\la(\mu)}(y_{1},\dotsc,y_{m}) = s_{\la(\mu)}(y_{1},\dotsc,y_{m})$.

Consider the restriction of Mason's analogue of
RSK to lower triangular matrices $A$.
In this case we have $m=n$ and the
following identity, studied by Lascoux in~\cite{Lascoux03}:
\begin{equation*}
\sum_{\mu\in \ZZ^n_{\ge0}} \bar{\kappa}^{\mu}(x_{1},\dotsc,x_{n}) \kappa_\mu(y_{1},\dotsc,y_{n})  = \prod_{1\le j\le i\le n}  \frac1{1-x_i y_j}.  
\end{equation*}
Define a probability measure on the space of all compositions $\mu$:
\begin{equation*}
\PP_n(\mu) = \prod_{1\le j\le i\le n}  (1-x_i y_j) \times \bar{\kappa}^{\mu}(x_{1},\dotsc,x_{n}) \kappa_\mu(y_{1},\dotsc,y_{n} ),\quad \mu\in \ZZ^n_{\ge0}.  
\end{equation*}

\begin{figure}[!h]
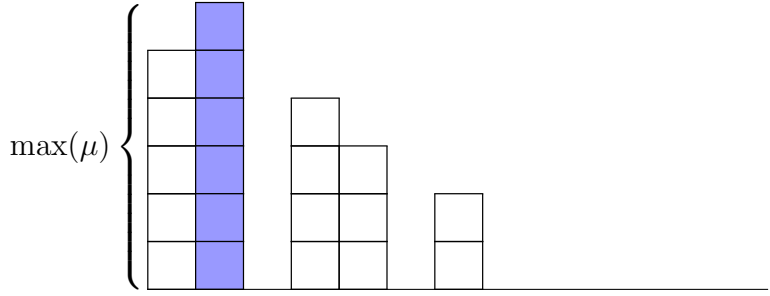

    \centering
    
    \begin{tabular}{r@{}l}
    \raisebox{-7.7ex}{$\max(\mu)\left\{\vphantom{\begin{array}{c}~\\[15.5ex] ~
    \end{array}}\right.$} &
    
    \begin{ytableau}
    \none & *(blue!40)  & \none  & \none  & \none  & \none \\
          & *(blue!40)  & \none  & \none  & \none  & \none \\
          & *(blue!40)  & \none  &        & \none  & \none \\
          & *(blue!40)  & \none  &        &        & \none \\
          & *(blue!40)  & \none  &        &        & \none & \; \\
          & *(blue!40)  & \none  &        &        & \none & \; \\        
        \end{ytableau}\\[-1.9ex]
      &\rule{0.5\linewidth}{0.4pt}
    \end{tabular}
    \caption{A composition $\mu\in\ZZ^{n}$, top entry is shown in blue}
    \label{fig:composition+top-entry}
\end{figure}

A natural question is to describe the asymptotic behavior of the top entry of such a composition,
\begin{equation*}
\PP(L_n=k) = \prod_{1\le j\le i\le n}  (1-x_i y_j) \times \sum_{\mu\in \ZZ^n_{\ge0}, \max(\mu) = k} \bar{\kappa}^{\mu}(x_{1},\dotsc,x_{n}) \kappa_\mu(y_{1},\dotsc,y_{n} ).  
\end{equation*}
A composition $\mu\in\ZZ^{n}$ with the marked top entry (highest column) is shown in Figure~\ref{fig:composition+top-entry}.

The usual way to write RSK algorithm is to represent a lower triangular matrix $A=(a_{ij})$ as a collection of pairs $\binom{i}{j}$ appearing $a_{ij}$ times, ordered lexicographically. Then Schensted row insertion applied to this sequence $\begin{pmatrix} i_{1},i_{2},\dotsc\\j_{1},j_{2},\dotsc
\end{pmatrix}$ produces a pair of semistandard Young tableaux of the same shape, and the length of the first row is equal to the length of a longest non-decreasing subsequence in $(j_{1},j_{2},\dotsc)$ and the LPP time $L^{lt}_{(1,1)\to(n,n)}$. Mason's RSK uses the column insertion instead and produces the LPP time $L^{lt}_{(1,n)\to(n,1)}$. According to Theorem~4.2 from \cite{AGL} we have $\PP(L^{lt}_{(1,n)\to(n,1)}=k) = \PP(L_n=k)$.

\begin{figure}[!h]
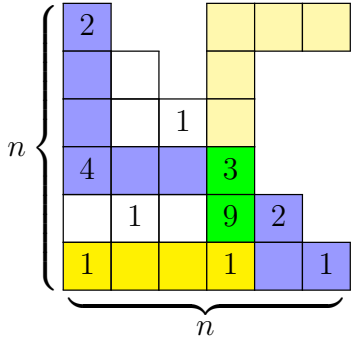

    \centering
    
    \begin{tabular}{r@{}l}
    \raisebox{-7.7ex}{$n\left\{\vphantom{\begin{array}{c}~\\[15.5ex] ~
    \end{array}}\right.$} &
    
    \begin{ytableau}
    *(blue!40) 2    & \none & \none & *(yellow!40)  & *(yellow!40) & *(yellow!40) \\
    *(blue!40)      & \; & \none      & *(yellow!40)      & \none \\
    *(blue!40)      & \; & 1 & *(yellow!40)      & \none \\
    *(blue!40) 4    & *(blue!40) & *(blue!40) & *(green) 3 & \none \\
          &     1    &            & *(green) 9 & *(blue!40) 2  & \none \\
     *(yellow) 1    &  *(yellow) &  *(yellow) & *(yellow) 1  & *(blue!40)  &     *(blue!40) 1   & \none \\        
    \end{ytableau}\\[-1.5ex]
    &\hspace{0.2em}$\underbrace{\hspace{8.8em}}_{\displaystyle n}$
    \end{tabular}
    \caption{Path from $(1,1)$ to $(n,n)$ is shown in blue, path from $(1,n)$ to $(n,1)$ in yellow, their intersection is shown in green. The continuation of a yellow path above the main diagonal is arbitrary and is shown in pale yellow. }
    \label{fig:two-paths}
\end{figure}

\begin{figure}[!h]
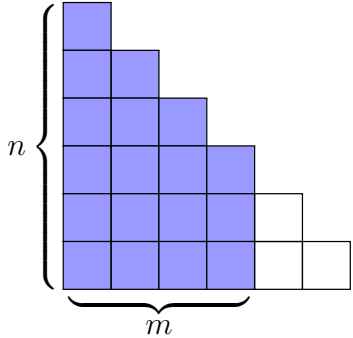

    \centering
    
    \begin{tabular}{r@{}l}
    \raisebox{-7.7ex}{$n\left\{\vphantom{\begin{array}{c}~\\[15.5ex] ~
    \end{array}}\right.$} &
    
    \begin{ytableau}
    *(blue!40)      & \none & \none & \none  & \none \\
    *(blue!40)      & *(blue!40) & \none      & \none      & \none \\
    *(blue!40)      & *(blue!40) & *(blue!40) & \none      & \none \\
    *(blue!40)      & *(blue!40) & *(blue!40) & *(blue!40) & \none \\
    *(blue!40)      & *(blue!40) & *(blue!40) & *(blue!40) &         & \none \\
    *(blue!40)      & *(blue!40) & *(blue!40) & *(blue!40) &         &        & \none \\        
    \end{ytableau}\\[-1.5ex]
    &\hspace{0.2em}$\underbrace{\hspace{5.8em}}_{\displaystyle m}$
    \end{tabular}
    \caption{Diagram of $\Lambda_{n, m}$}
    \label{truncation-shape}
\end{figure}

We study the asymptotic behavior for both LPP times  $L^{lt}_{(1,1)\to(n,n)}$ and  $L^{lt}_{(1,n)\to(n,1)}$ for random lower triangular matrices $A$. In Figure~\ref{fig:two-paths} we show both paths from $(1,1)$ to $(n,n)$ and from $(1,n)$ to $(n,1)$. 
In the former case we can also consider truncated staircase shape $\Lambda_{n,m}$ (see Fig.~\ref{truncation-shape}) and corresponding truncated lower-triangular matrices $A$ with non-zero elements appearing only inside $\Lambda_{n,m}$.
For the constant specialization $x_{i}=y_{j}=\sqrt{q}$ for all $i,j$ and $q\in(0,1)$ we can present a distribution of $L^{lt}_{(1,1)\to(n,m)}$ as a marginal of a Pfaffian Schur process. Such distributions can be analyzed by methods from \cite{OR03,BR01,BBNV}.
Moreover, the methods of \cite{BR01} can be applied to $L^{lt}_{(1,n)\to(n,1)}$ as well. 

We obtain the following results. Let $F_{GOE}(s)$ be the GOE Tracy--Widom distribution defined in~\eqref{eq:F^GOE} below and $F_{GSE}(s)$ be the GSE Tracy--Widom distribution defined in~\eqref{eq:F^GSE} below. Set $b_q = \frac{\sqrt{q}}{1-\sqrt{q}}$, 
$c_q = \frac{1-\sqrt{q}}{q^{1/6}(1+\sqrt{q})^{1/3}}$.

\begin{thm}\label{thm:asympt-staircase_Demazure}
    Consider the parameters $x_i=y_i=\sqrt{q}$, $i\ge1$, $q\in (0,1)$. We have
    \begin{equation*}
    \lim_{n\to\infty} \PP\left(L_n \le b_q n + 
        \frac {c_q^{-1}}{2} n^{1/3} s\right) = F_{GOE}(s).       
    \end{equation*}
\end{thm}

From~\cite{BR01} we have the following straightforward corollary.
\begin{prop}\label{prop:asympt-staircase_Pfaffian}
    Consider the parameters $x_i=y_i=\sqrt{q}$, $i\ge1$, $q\in (0,1)$. We have
    \begin{equation*}
    \lim_{n\to\infty} \PP(L^{lt}_{(1,1)\to(n,n)} \le 2 b_q n + c_q^{-1} n^{1/3} s) = F_{GSE}(s).       
    \end{equation*}

\end{prop}

The same limiting behavior holds for $L^{lt}_{(1,1)\to(n,m)}$ if $n-m = o(n^{2/3})$. A non-trivial change in asymptotics appears for  $n-m = \bigO(n^{2/3})$. Let $F_{u,\infty}(s)$ be a generalisation of GSE Tracy--Widom distribution defined in~\eqref{eq:gen-F_GSE} below.

\begin{figure}[h]
  \includegraphics[width=0.45\linewidth]{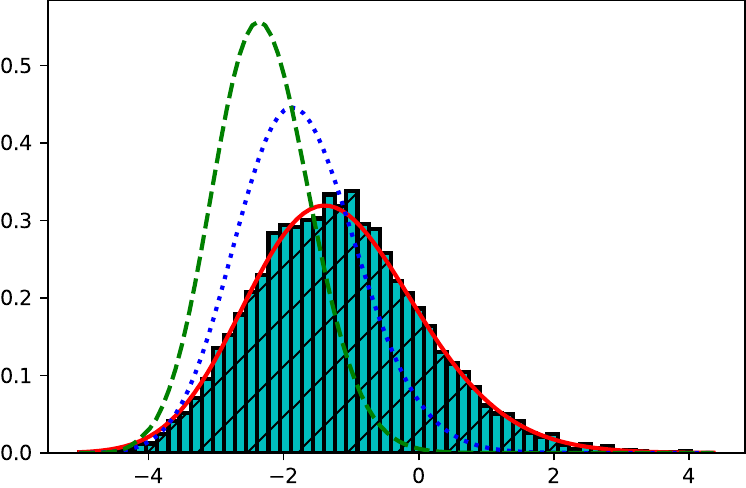}
  \includegraphics[width=0.45\linewidth]{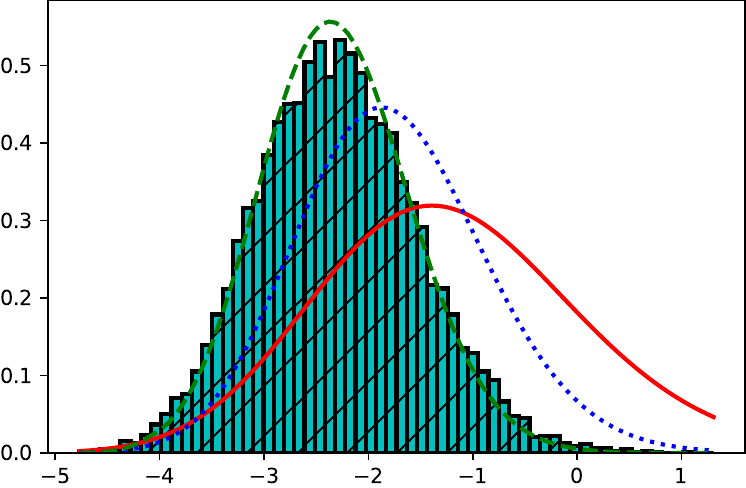}    
  \caption{{\it Left:} Histogram of LPP time $L^{lt}_{(1,n)\to(n,1)}$ with $10000$ samples for $n=100$ and $q=0.618$ and GSE, GUE, GOE Tracy-Widom distributions in green, blue and red correspondingly.\\
    {\it Right:} Histogram of LPP time $L^{lt}_{(1,1)\to(n,n)}$ with $10000$ samples for $n=100$ and $q=0.618$ and GSE, GUE, GOE Tracy-Widom distributions in green, blue and red correspondingly.}
  \label{fig:histograms}
\end{figure}

\begin{thm}\label{thm:asympt-trunc_staircase}
  Consider the parameters $x_i=y_i=\sqrt{q}$, $i\ge1$, $q\in (0,1)$. For $u\ge 0$ we have
  \begin{equation*}
    \lim_{n\to\infty} 
    \PP(L^{lt}_{(1,1)\to(n,n-u n^{2/3})} \le 
    b_q (2n - u n^{2/3})  + c_q^{-1} n^{1/3} s) 
    = F_{u, \infty}(s).      
  \end{equation*}
\end{thm}
In Figure~\ref{fig:histograms} we present simulation results for LPP times  $L^{lt}_{(1,n)\to(n,1)}$ and  $L^{lt}_{(1,1)\to(n,n)}$ obtained by applying RSK algorithm to random lower triangular matrices with $10000$ samples for $n=100$ and $q=0.618$, and graphs for the densities of GSE, GUE, GOE Tracy-Widom distributions for comparison.

\section*{Acknowledgments}
\label{sec:acknowledgements}
We thank Ievgen Makedonskyi, Nicolai Reshetikhin and Travis Scrimshaw for fruitful conversations.

\section{Proofs}

According to Theorem 4.2 from~\cite{AGL}, we have
 $ \PP(L^{lt}_{(1,n)\to(n,1)}=k) = \PP(L_n=k). $

Now consider a symmetric random $n\times n$ matrix $A = (a_{ij})$ with entries given by $a_{ij} = a_{ji} = w_{i,j}$, where $w_{i,j}$, $j\le i$ are the geometric random variables from ~\eqref{eq:square-LPP-def}. Define $L^{sym}_{(1,n)\to(n,1)}$ to be the symmetric LPP time, i.e., 
\begin{equation*}
L^{sym}_{(1,n)\to(n,1)} = \max_{\pi}\sum_{(r,t)\in\pi}a_{r,t},  
\end{equation*}
where we maximize over all the paths from $(1,n)$ to $(n,1)$.
Using the symmetry with respect to the main diagonal, we immediately see that the distributions of $L^{lt}_{(1,n)\to(n,1)}$ and $\frac12 L^{sym}_{(1,n)\to(n,1)}$ are the same.

Consider a similar last passage percolation (LPP) model. Take $(a_n)_{n\in\NN}$, $a_n\in (0,1)$, $n\in\NN$, 
$\alpha\in \left(0,\inf_{n\in\NN}(a_n^{-1})\right)$. Let $\{\tilde w_{i,j}\}_{1\le i,j\le n}$ be a family of independent geometric random variables,
\begin{equation}\label{eq:LPP_alpha}
  \begin{aligned}
    &\PP(\tilde w_{i,j} = k) = (1-a_i a_j) (a_i a_j)^k,  & k\in \ZZ_{\ge 0}, &\quad i\neq j,  \\
    &\PP(\tilde w_{i,j} = k) = (1-\alpha a_i) (\alpha a_i)^k,  & k\in \ZZ_{\ge 0},  &\quad i=j.    
  \end{aligned}
\end{equation}
Denote by $L^{sym, \alpha}_{(1,n)\to(n,1)}$ the corresponding 
symmetric LPP time from $(1,n)$ to $(n,1)$.

Asymptotics for $L^{sym, \alpha}_{(1,n)\to(n,1)}$ for the constant specialization $a_i=\sqrt{q}$, $q\in (0,1)$, $i\ge1$, is given by the GOE Tracy--Widom distribution ~\cite[Theorem 4.2]{BR}:
\begin{equation*}
  \lim_{n\to\infty} \PP(L^{sym, \alpha}_{(1,n)\to(n,1)} \le 2 b_q n + c_q^{-1} n^{1/3} s) = F_{GOE}(s),
\end{equation*}
where $F_{GOE}$ is defined as follows. 

Let $q(x)$ be the solution of the Painlev\'e II equation 
$$
q''(x) = 2q^2(x) + xq(x),
$$
with the asymptotic behavior $q(x) = (1+o(1)) \Ai(x)$, $x\to\infty$, where $\Ai(x)$ is the Airy function. Define the Tracy--Widom distribution as
$$
F_{GUE}(s) = \exp(-\int_s^\infty (x-s)q^2(x)dx),
$$
the GOE Tracy--Widom distribution as
\begin{equation}\label{eq:F^GOE}
   F_{GOE}(s) = (F_{GUE}(s))^{1/2} \exp\left( - \frac12 \int_s^\infty q(x)dx \right), 
\end{equation}
and the GSE Tracy--Widom distribution as
\begin{equation}\label{eq:F^GSE}
   F_{GSE}(s) = (F_{GUE}(s))^{1/2} \cosh\left(\frac12 \int_s^\infty q(x)dx \right).
\end{equation}

%Let $\mcC_a^{\phi}$ be the contour, defined as a union of two semi-infinite rays from $a\in\CC$ with angles $\phi$ and $-\phi$, oriented from $a+\infty e^{-i\phi}$ to $a+\infty e^{i\phi}$. Define 
%\begin{align*}
%    K^{GOE}_{11}(\xi, \xi') &= 
%    \int_{\mcC_1^{\pi/3}} \frac{dz}{2\pi i}
%    \int_{\mcC_1^{\pi/3}} \frac{dw}{2\pi i}
%    \frac{z-w}{z+w}e^{z^3/3+w^3/3-z\xi-w\xi'}, \\
%     K^{GOE}_{12}(\xi, \xi') &= 
%     - K^{GOE}_{21}(\xi, \xi') =
%    \int_{\mcC_1^{\pi/3}} \frac{dz}{2\pi i}
%    \int_{\mcC_{-1/2}^{\pi/3}} \frac{dw}{2\pi i}
%    \frac{w-z}{2w(z+w)}e^{z^3/3+w^3/3-z\xi-w\xi'},\\
%     K^{GOE}_{22}(\xi, \xi') &= 
%    \int_{\mcC_1^{\pi/3}} \frac{dz}{2\pi i}
%    \int_{\mcC_1^{\pi/3}} \frac{dw}{2\pi i}
%    \frac{z-w}{4zw(z+w)}e^{z^3/3+w^3/3-z\xi-w\xi'} \\
%    &+ \int_{\mcC_1^{\pi/3}} \frac{dz}{8z\pi i} e^{z^3/3-z\xi}
%    - \int_{\mcC_1^{\pi/3}} \frac{dz}{8z\pi i} e^{z^3/3-z\xi'}
%    - \frac{1}{4}\II_{x>y} + \frac{1}{4} \II_{y<x}.
%\end{align*}

%and
%\begin{equation}\label{eq:F^GOE}
%    F_{GOE}(s) = \Pf(J+K^{GOE})_{(s,\infty)}, \quad
%    K^{GOE}(\xi, \xi') = 
%\begin{pmatrix}
%    K^{GOE}_{11}(\xi, \xi')  & K^{GOE}_{12}(\xi, \xi') \\
%    -K^{GOE}_{12}(\xi, \xi') & K^{GOE}_{22}(\xi, \xi')    
%\end{pmatrix}.
%\end{equation}

Now consider $L^{sym, \alpha}_{(1,n)\to(n,1)}$ for the constant specialization $a_i=\sqrt{q}$, $i\ge1$, $q\in (0,1)$, $\alpha = \sqrt{q}$. We see that this distribution is equal to $L^{sym}_{(1,n)\to(n,1)}$ for the constant specialization $x_i=y_i=\sqrt{q}$, and Theorem~\ref{thm:asympt-staircase_Demazure} is proved.

We now proceed to the proof of Proposition~\ref{prop:asympt-staircase_Pfaffian} and Theorem~\ref{thm:asympt-trunc_staircase}. Now for the same lower triangular LPP model consider the symmetric LPP time from $(1,1)$ to $(n,n)$: 
$$
L^{sym}_{(1,1)\to(n,n)} = \max_{\pi}\sum_{(r,t)\in\pi}a_{r,t}.
$$
Using the symmetry with respect to the main diagonal, we immediately see that the distributions of $L^{lt}_{(1,1)\to(n,n)}$ and $L^{sym}_{(1,1)\to(n,n)}$ are the same.

For last passage percolation (LPP) model~\eqref{eq:LPP_alpha}
denote by $L^{sym, \alpha}_{(1,1)\to(n,n)}$ the corresponding symmetric LPP time from $(1,1)$ to $(n,n)$.

Asymptotics for $L^{sym, \alpha}_{(1,1)\to(n,n)}$ for the constant
specialization $a_i=\sqrt{q}$, $i\ge1$, $q\in (0,1)$, $\alpha<1$ is
given by the GSE Tracy--Widom distribution ~\cite[Theorem 4.2]{BR}:
\begin{equation*}
  \lim_{n\to\infty} \PP(L^{sym, \alpha}_{(1,1)\to(n,n)} \le 2 b_q n + c_q^{-1} n^{1/3} s) = F_{GSE}(s),
\end{equation*}
where $F_{GSE}$ is defined in~\eqref{eq:F^GSE}.  
%Define 
%\begin{equation*}
%    K^{GSE}_{11}(\xi, \xi') = 
%    \int_{\mcC_1^{\pi/3}} \frac{dz}{2\pi i}
%    \int_{\mcC_1^{\pi/3}} \frac{dw}{2\pi i}
%    \frac{z-w}{4zw(z+w)}e^{z^3/3+w^3/3-z\xi-w\xi'},
%  \end{equation*}
%  \begin{equation*}
%     K^{GSE}_{12}(\xi, \xi') = 
%     - K^{GSE}_{21}(\xi, \xi') =
%    \int_{\mcC_1^{\pi/3}} \frac{dz}{2\pi i}
%    \int_{\mcC_1^{\pi/3}} \frac{dw}{2\pi i}
%    \frac{z-w}{4z(z+w)}e^{z^3/3+w^3/3-z\xi-w\xi'},
%  \end{equation*}
%  \begin{equation*}
%     K^{GSE}_{22}(\xi, \xi') = 
%    \int_{\mcC_1^{\pi/3}} \frac{dz}{2\pi i}
%    \int_{\mcC_1^{\pi/3}} \frac{dw}{2\pi i}
%    \frac{z-w}{4(z+w)}e^{z^3/3+w^3/3-z\xi-w\xi'},  
%\end{equation*}
%and
%\begin{equation}\label{eq:F^GSE}
%    F_{GSE}(s) = \Pf(J+K^{GSE})_{(s,\infty)}, \quad
%    K^{GSE}(\xi, \xi') = 
%\begin{pmatrix}
%    K^{GSE}_{11}(\xi, \xi')  & K^{GSE}_{12}(\xi, \xi') \\
%    -K^{GSE}_{12}(\xi, \xi') & K^{GSE}_{22}(\xi, \xi')    
%\end{pmatrix}.
%\end{equation}

Now consider $L^{sym, \alpha}_{(1,1)\to(n,n)}$ for $\alpha = \sqrt{q}$. We see that this distribution is equal to $L^{sym}_{(1,1)\to(n,n)}$ for the constant specialization 
$x_i=y_i=\sqrt{q}$. Thus Proposition~\ref{prop:asympt-staircase_Pfaffian} is proved.

Finally we consider the paths from $(1,1)$ to $(n,m)$. Once again, using the symmetry with respect to the main diagonal, we see that the distributions of $L^{lt}_{(1,1)\to(n,m)}$ and $L^{sym}_{(1,1)\to(n,m)}$ are the same. Substituting constant specialization $a_i=\sqrt{q}$, $\alpha=\sqrt{q}$ into Theorem~4.1 from~\cite{BBNV} we obtain the following result. Set
$$
J(x,y) = \II_{x=y} 
\begin{pmatrix}
    0  & 1 \\
    -1 & 0
\end{pmatrix}.
$$
\begin{prop}\label{thm:L-nm_distr}
Consider the parameters $x_i=y_i=\sqrt{q}$, $i\ge1$, $q\in (0,1)$. For $m,n\in\NN$, $m\le n$ we have
$$
\PP\left( L^{lt}_{(1,1)\to(n,m)} \le k \right) =
    \Pf(J - K)_{(k,\infty)},
$$
where the kernel $K$ is given by
$$
    K(x_1, x_2) = 
\begin{pmatrix}
    K_{11}(x_1,x_2)  & K_{12}(x_1,x_2) \\
    -K_{12}(x_1,x_2) & K_{22}(x_1,x_2)    
\end{pmatrix},
$$
\begin{multline*}
K_{11}(x_1,x_2) = 
    \oint \frac{dz}{2\pi i} \oint \frac{dw}{2\pi i} 
    \frac{z-w}{(z^2-1)(w^2-1)(zw-1)}z^{-x_1+1/2}w^{-x_2+1/2} \\
    \times \prod_{i=1}^{n}\frac1{1-\sqrt{q}z}
    \prod_{i=1}^{n}\frac1{1-\sqrt{q}w}
    \prod_{i=1}^{m+1}(1-\sqrt{q}/z)
    \prod_{i=1}^{m+1}(1-\sqrt{q}/w),
\end{multline*}
for counterclockwise oriented circle contours around $0$ satisfying $1<|z|$, $|w|<\frac1{\sqrt{q}}$;
\begin{multline*}
K_{12}(x_1,x_2) = 
    \oint \frac{dz}{2\pi i} \oint \frac{dw}{2\pi i} 
    \frac{zw-1}{(z-w)(z^2-1)}z^{-x_1+1/2}w^{x_2-3/2} \\
    \times \prod_{i=1}^{n}\frac1{1-\sqrt{q}z}
    \prod_{i=1}^{n}(1-\sqrt{q}w)
    \prod_{i=1}^{m+1}(1-\sqrt{q}/z)
    \prod_{i=1}^{m+1}\frac1{1-\sqrt{q}/w},
\end{multline*}
for counterclockwise oriented circle contours around $0$ satisfying $\sqrt{q}<|w|<|z|<\frac1{\sqrt{q}}$ and $1<|z|$; and 
\begin{multline*}
K_{22}(x_1,x_2) = 
    \oint \frac{dz}{2\pi i} \oint \frac{dw}{2\pi i} 
    \frac{z-w}{zw-1}z^{x_1-3/2}w^{x_2-3/2} \\
    \times \prod_{i=1}^{n}(1-\sqrt{q}z)
    \prod_{i=1}^{n}(1-\sqrt{q}w)
    \prod_{i=1}^{m+1}\frac1{1-\sqrt{q}/z}
    \prod_{i=1}^{m+1}\frac1{1-\sqrt{q}/w},
\end{multline*}
for counterclockwise oriented circle contours around $0$ satisfying $\sqrt{q}<|z|$, $|w|<\frac1{\sqrt{q}}$ and $1<|zw|$.
\end{prop}

To the best of our knowledge, the asymptotic behaviour for $L^{sym, \alpha}_{(1,1)\to(n,m)}$, $m\neq n$, for constant specialization, with $\alpha < 1$, haven't been studied. But there is similar asymptotic result in~\cite{BBNV}, Theorem~4.2, dealing with the case $\alpha = 1 - 2vc_qn^{-1/3}$, $v\in\RR$, that we state below. 

Let $\mcC_a^{\phi}$ be the contour, defined as a union of two semi-infinite rays from $a\in\CC$ with angles $\phi$ and $-\phi$, oriented from $a+\infty e^{-i\phi}$ to $a+\infty e^{i\phi}$.

\begin{thm}{\cite[Theorem 4.2]{BBNV}}
    Consider the LPP time $L^{sym, \alpha}_{(1,1)\to(n,m)}$ for the model~\eqref{eq:LPP_alpha} with the parameters $a_i=\sqrt{q}$, $i\ge1$, $q\in (0,1)$, $\alpha = 1 - 2vc_qn^{-1/3}$, $v\in\RR$, $m=n-u n^{2/3}$. For $u > 0$ we have
    \begin{equation*}
    \lim_{n\to\infty} \PP(L^{sym, \alpha}_{(1,1)\to(n,n-u n^{2/3})} \le 
    b_q (2n - u n^{2/3})  + c_q^{-1} n^{1/3} s) 
    = F_{u, v}(s), 
    \end{equation*}
where 
\begin{equation*}
    F_{u, v}(s) = \Pf(J+K^{u,v})_{(s,\infty)}, \quad
    K^{u,v}(\xi, \xi') = 
\begin{pmatrix}
    K^{u,v}_{11}(\xi, \xi')  & K^{u,v}_{12}(\xi, \xi') \\
    -K^{u,v}_{12}(\xi, \xi') & K^{u,v}_{22}(\xi, \xi')    
\end{pmatrix}.
\end{equation*}
  and

  \begin{equation*}
    K^{u,v}_{11}(\xi, \xi') = 
    \int_{\mcC_1^{\pi/3}} \frac{dz}{2\pi i}
    \int_{\mcC_1^{\pi/3}} \frac{dw}{2\pi i}
    \frac{(z-w)(w+2v)(z+2v)}{4zw(z+w)}e^{z^3/3+w^3/3-u(z^2+w^2)-z\xi-w\xi'},
  \end{equation*}
  \begin{equation*}
     K^{u,v}_{12}(\xi, \xi') = 
     - K^{u,v}_{21}(\xi, \xi') =
    \int_{\mcC_{A_z}^{\pi/3}} \frac{dz}{2\pi i}
    \int_{\mcC_{A_w}^{2\pi/3}} \frac{dw}{2\pi i}
    \frac{(z+w)(z+2v)}{2(w+2v)z(z-w)}e^{z^3/3-w^3/3-u(z^2-w^2)-z\xi+w\xi'},
  \end{equation*}
  \begin{equation*}
     K^{u,v}_{22}(\xi, \xi') = K^{u,v,1}_{22}(\xi, \xi') + K^{u,v,2}_{22}(\xi, \xi') 
  \end{equation*}
  \begin{equation*}
     K^{u,v,1}_{22}(\xi, \xi') = 
    \int_{\mcC_{B_3}^{2\pi/3}} \frac{dz}{2\pi i}
    \int_{\mcC_{B_4}^{2\pi/3}} \frac{dw}{2\pi i}
    \frac{z-w}{(w+2v)(z+2v)(z+w)}e^{-z^3/3-w^3/3+u(z^2+w^2)+z\xi+w\xi'},  
  \end{equation*}
  \begin{multline*}
    K^{u,v,2}_{22}(\xi, \xi') =
    e^{8v^3/3 +4v^2u-2v\xi'}
    \int_{\mcC_{B_2}^{2\pi/3}} \frac{dz}{2\pi i}
    \frac{1}{z-2v}e^{-z^3/3+uz^2+z\xi}   -  
    \\
    \int_{\mcC_{B_1}^{2\pi/3}} \frac{dz}{2\pi i}
    \frac{2z}{(z+2v)(z-2v)}e^{2uz^2+z(\xi-\xi')},    
  \end{multline*}
where $A_z> A_w>-2v$, $A_z>0$; and $B_1>2|v|$, $B_2>2v$, $B_3>-2v>B_4$, $B_3<-B_4$. 
\end{thm}

We can obtain similar result in our case. For $u\in\RR$ we define 
\begin{align*}
    K^{u,\infty}_{11}(\xi, \xi') = 
    \int_{\mcC_1^{\pi/3}} \frac{dz}{2\pi i}
    \int_{\mcC_1^{\pi/3}} \frac{dw}{2\pi i}
    \frac{z-w}{4zw(z+w)}e^{z^3/3+w^3/3-u(z^2+w^2)-z\xi-w\xi'}, \\
     K^{u,\infty}_{12}(\xi, \xi') = 
     - K^{u,\infty}_{21}(\xi, \xi') =
    \int_{\mcC_1^{\pi/3}} \frac{dz}{2\pi i}
    \int_{\mcC_1^{\pi/3}} \frac{dw}{2\pi i}
    \frac{z-w}{2z(z+w)}e^{z^3/3+w^3/3-u(z^2-w^2)-z\xi-w\xi'},\\
     K^{u,\infty}_{22}(\xi, \xi') = 
    \int_{\mcC_1^{\pi/3}} \frac{dz}{2\pi i}
    \int_{\mcC_1^{\pi/3}} \frac{dw}{2\pi i}
    \frac{z-w}{(z+w)}e^{z^3/3+w^3/3+u(z^2+w^2)-z\xi-w\xi'},  
\end{align*}
and set
\begin{equation}\label{eq:gen-F_GSE}
    F_{u, \infty}(s) = \Pf(J+K^{u,\infty})_{(s,\infty)}, \quad
    K^{u,\infty}(\xi, \xi') = 
\begin{pmatrix}
    K^{u,\infty}_{11}(\xi, \xi')  & K^{u,\infty}_{12}(\xi, \xi') \\
    -K^{u,\infty}_{12}(\xi, \xi') & K^{u,\infty}_{22}(\xi, \xi')    
\end{pmatrix}.
\end{equation}

One can define $F_{u,v}(s)$ for $u=0$, see~\cite{BBNV} for details. The GSE Tracy--Widom distribution admits a similar Pfaffian representation~\cite{BBCS2018}. And we can obtain GSE Tracy--Widom distribution $F_{GSE}(s)$ as a limit of $F_{0, v}(s)$, $v\to\infty$. Similarly, one can think about the function $F_{u, \infty}(s)$ as a $v=\infty$ case of the function $F_{u,v}(s)$. Asymptotic analysis similar to Theorem~4.2 from~\cite{BBNV} and \cite{BBCS2018} leads to the Theorem~\ref{thm:asympt-trunc_staircase}. We plan to present a detailed proof of this theorem in a subsequent paper.

\section{Conclusion and open problems}
\label{sec:concl-open-probl}

In Theorem~\ref{thm:asympt-staircase_Demazure} we have established the asymptotics for the top entry of a random composition with respect to a Demazure measure. As far as we know local and global asymptotic behavior of Demazure measures was not considered previously. It is possible to consider a Demazure measure on diagrams by ordering elements of a random composition in non-increasing order. We have numerical evidence for convergence of such random diagrams to a limit shape. For constant specialization there is a natural description of the measure in representation-theoretic terms, so it could be possible to obtain a closed formula for the measure itself and solve a variational problem for the limit shape.

Similar non-symmetric Cauchy identities for a general shape $\Lambda \subset (n^n)$ were introduced by Feigin, Khoroshkin and Makedonskyi in 
\cite{feigin2024cauchyidentitiesstaircasematrices,khoroshkin2025bubblesorthoweduality}. We expect similar limit shape convergence to hold in this case also.

Proposition~\ref{prop:asympt-staircase_Pfaffian} and Theorem~\ref{thm:asympt-trunc_staircase} use the Pfaffian structure of the model to obtain the asymptotic distribution of the LPP time, or the length of the first row of a random Young diagram. Such a structure appears because of the correspondence between the lower-triangular and the symmetric models. This correspondence breaks if we consider the second and the following rows. Therefore local asymptotics of the lower-triangular model remains an open problem.

Demazure measures are based on the column insertion while Proposition~\ref{prop:asympt-staircase_Pfaffian} and Theorem~\ref{thm:asympt-trunc_staircase} are based on the row insertion. It is natural to expect that these results might be useful for the study of the asymptotic behavior of a dual version of a Demazure measure.

\bibliographystyle{alpha}
\bibliography{shapes}{}

\end{document}